\documentclass[12pt]{article}
\usepackage{amsmath,amssymb,amsbsy,amsfonts,amsthm,
            latexsym,amsopn,amstext,amsxtra,euscript,amscd}

\begin{document}

\newtheorem{prop}{Proposition}
\newtheorem{lemma}[prop]{Lemma}
\newtheorem{cor}[prop]{Corollary}
\newtheorem{conj}[prop]{Conjecture}
\newtheorem{defi}[prop]{Definition}
\newtheorem{theorem}[prop]{Theorem}
\newtheorem{fac}[prop]{Fact}
\newtheorem{facs}[prop]{Facts}
\newtheorem{com}[prop]{Comments}
\newtheorem{prob}{Problem}
\newtheorem{problem}[prob]{Problem}
\newtheorem{ques}{Question}
\newtheorem{question}[ques]{Question}
\newtheorem{remark}[prop]{Remark}

\def\em{{\mathbf{e}}_m}
\def\ek{{\mathbf{e}}_k}
\def\ep{{\mathbf{e}}_p}
\def\et{{\mathbf{e}}_{t_p}}
\def\etd{{\mathbf{e}}_{t_p/d}}
\def\erd{{\mathbf{e}}_{r/d}}
\def\cS{\mathcal{S}}
\def\cE{\mathcal{E}}
\def\cX{\mathcal{X}}
\def\cP{\mathcal{P}}
\def\cK{\mathcal{K}}
\def\eT{{\mathbf{e}}_T}

\title{\bf The large sieve for $2^{[O(n^{15/14+o(1)})]}$ modulo primes}
\author{
{\sc M.~Z.~Garaev}\\
{\normalsize Instituto\ de Matem{\'a}ticas,
UNAM} \\
{\normalsize Campus Morelia, Ap. Postal 61-3 (Xangari)}\\
{\normalsize C.P.\ 58089, Morelia, Michoac{\'a}n, M{\'e}xico} \\
{\normalsize \tt garaev@matmor.unam.mx}}

\date{}

\maketitle

\begin{abstract}
Let $\lambda$ be a fixed integer, $\lambda\ge 2.$ Let $s_n$ be any
strictly increasing sequence of positive integers satisfying $s_n\le
n^{15/14+o(1)}.$ In this paper we give a version of the large sieve
inequality for the sequence $\lambda^{s_n}.$ In particular, we prove
that for $\pi(X)(1+o(1))$ primes $p, \ p\le X,$ the numbers
$$
\lambda^{s_n},\quad n\le X(\log X)^{2+\varepsilon}
$$
are uniformly distributed modulo $p.$
\end{abstract}



\newpage

\section{Notation}

Throughout the paper the following notations will be used:

\bigskip

$\lambda$ denotes a fixed positive integer, $\lambda\ge 2$;

\bigskip

$X$ and $T$ are large parameters, $T$ is an integer;

\bigskip

$\Delta>X^{1/3}$ is a parameter;

\bigskip

$s_n, \ n=1,2,\ldots,$ is a strictly increasing sequence of positive
integers (which may depend on the parameters $X, T, \Delta$);

\bigskip

$\gamma_n, \, n=1,2,\ldots,$ are any complex coefficients (which may
depend on the parameters $X, T, \Delta$) with $|\gamma_n|\le 1$;

\bigskip

$p$ and $q$ always denote prime numbers;

\bigskip

$t_p$  denotes  the multiplicative order of $\lambda$ modulo $p$;

\bigskip

$ \cE=\cE(\Delta, X)=\{p: p\le X, t_p>\Delta\};$ that is the set of
all primes $p,\ p \le X,$ with  $ t_p> \Delta;$

\bigskip

 For integers $a$ and $b$, their
greatest common divisor is denoted by $(a,b).$

Given a set $\cX$ we use $|\cX|$ to denote its cardinality.

As usual, $\pi(X)$ denotes the number of primes not exceeding $X,$
and $\tau(n)$ denotes the number of positive integer divisors of
$n$. We also follow the standard abbreviation
$$
\em(z)=e^{2\pi iz/m}.
$$

\section{Introduction}

Recently, J. Bourgain~\cite{Bourg, Bourg1} has proved that for
$\pi(X)(1+o(1))$ primes $p, p\le X,$ the Mersenne numbers
$M_q=2^q-1,$ $q\le X^{2+\varepsilon},$ are uniformly distributed
modulo $p$ for any given $\varepsilon>0.$ Furthermore, he has
explicitly described the set of primes $p$ for which we can be sure
that the Mersenne numbers are uniformly distributed modulo $p.$ This
set is expressed in terms of certain conditions to the size of the
multiplicative order of $2$ modulo $p,$ which are satisfied for
almost all primes $p.$

Bourgain's result is based on his deep work related to nontrivial
estimates of double trigonometric sums. The possibility of
applications of such estimates to investigate Mersenne numbers in
residue classes modulo $p$ has been first discovered in~\cite{BCFS}.

An alternative approach, based on the large sieve inequality, has
been recently suggested in~\cite{GaSh}. From the result of Erd\H{o}s
and Murty~\cite{ErdMur} we know that the estimate $t_p>X^{1/2+o(1)}$
holds for almost all primes $p, p\le X.$ This has been used
in~\cite{GaSh} to obtain a nontrivial bound for the exponential sum
$$
\max_{(a, p)=1}\Big|\sum_{n\in \cS_N}e^{2\pi i
\frac{a\lambda^{n}}{p}}\Big|
$$
for $\pi(X)+o(\pi(X))$ primes $p, p\le X,$ provided that
$\cS_N\subset [1, N]$ is sufficiently dense (that is
$|\cS_N|>N^{1+o(1)}$) and $N$ is of the size $X^{1+o(1)}.$

The result of~\cite{GaSh} does not apply for sparser sets $\cS_N,$
but it is shown that such results can be obtained conditionally,
namely assuming the truth of the Extended Riemann Hypothesis.

In the present paper we provide a new argument which allows to deal
with sparse sets $\cS_N$ unconditionally. In particular, we obtain
equidistribution properties of $\lambda^n\pmod p,\ n\in \cS_N$ with
$|\cS_N|>N^{14/15-o(1)}.$ We show that further improvement could be
obtained if one knows how to complement in appropriate way the set
of exponent pairs for Gauss sums obtained by Konyagin.

Furthermore, while the result of~\cite{GaSh} only apply for the set
of primes $p\le X$ with $t_p>X^{1/2}(\log X)^c, \ c>0,$ here our
result works when $t_p>\Delta,$ where, depending on how sparse the
set $\cS$ is, $\Delta$ varies in $(X^{1/3+\varepsilon},
X^{1/2+o(1)}].$ This is useful if one is interested in obtaining
sharp upper bound estimates for the exceptional set of primes $p$ in
the equidistribution problem of the sequence $\lambda^n\pmod p, \,
n\in \cS_N.$

In what follows, we use the Landau symbol `$o$', as well as the
Vinogradov symbols `$\ll$' and `$\gg$' in their usual meanings. The
implied constants may depend on the small positive quantity
$\varepsilon,$ $\lambda$  and other fixed constants, and also on the
choice of the function $\nu(n)$  (in
Corollary~\ref{cor:Mersennetype} below, see also~\eqref{eqn:s_n}).

\section{Results}

The following statement is the main result of our paper. We recall
that $s_n,\ n=1,2\ldots,$ is any sequence of strictly increasing
positive integers.

\begin{theorem}
\label{thm:main1} For any $L>0$ the following bound holds:
\begin{eqnarray*}
\sum_{p\in \cE}\frac{1}{\tau(p-1)}\max_{(a, p)=1}\Big|\sum_{n\le
T}\gamma_n\ep(a\lambda^{s_n})\Big|^2 \ll
\left(X+s_TX^{1/7}\Delta^{-3/7}L+TL^{-7/4}\right)XT.
\end{eqnarray*}
\end{theorem}

If we optimize the choice of $L,$ then the estimate can be
reformulated in the form
$$
\sum_{p\in \cE}\frac{1}{\tau(p-1)}\max_{(a, p)=1}\Big|\sum_{n\le T}\gamma_n\ep(a\lambda^{s_n})\Big|^2\\
\ll \left(1+(s_T^7T^4X^{-10}\Delta^{-3})^{1/11}\right)X^2T.
$$

 As we have already mentioned in the Introduction, for
$\pi(X)(1+o(1))$ primes $p,\ p\le X,$ the inequality
$t_p>X^{1/2+g(X)}$ holds for any given function $g(x)=o(1).$

Let now $s_n$ satisfy the condition
\begin{equation}
\label{eqn:s_n} s_n\le n^{15/14+\nu_n}, \quad
\lim_{n\to\infty}\nu_n=0,
\end{equation}
where $\nu_n$ is an absolutely fixed sequence (therefore, does not
depend on the parameters $T, X, \Delta$). Set $T=[X(\log
X)^{2+\varepsilon}]$ and  take
$$
L=T^{|\nu_T|}(\log T)^{10},\quad \Delta=T^{1/2}L^7.
$$
Obviously, $L^7=X^{o(1)},\, \Delta=X^{1/2+o(1)}.$ Therefore,
$$
|\cE|=\pi(X)(1+o(1)).
$$
Incorporating this choice of the parameters in
Theorem~\ref{thm:main1}, we obtain
\begin{eqnarray*}
&&\sum_{p\in\cE}\frac{1}{\tau(p-1)}\max_{(a, p)=1}\Big|\sum_{n\le
T}\gamma_n\ep(a\lambda^{s_n})\Big|^2\\
&&\ll X^2T+XT^2(\log T)^{-10}\ll XT^2(\log T)^{-2-\varepsilon}.
\end{eqnarray*}

Next, let $\cE'$ be the subset of $\cE$ with $\tau(p-1)<(\log
X)^{1+\varepsilon/2}.$ From the Titchmarsh bound
\begin{equation}
\label{eqn:titch} \sum_{p\le X}\tau(p-1)\ll X
\end{equation}
(see for example Theorem~7.1 in Chapter~5 of~\cite{Prach}) it
follows that the inequality
$$
\tau(p-1)\ll (\log X)^{1+\varepsilon/2}
$$
holds for $\pi(X)(1+O((\log X)^{-\varepsilon/2}))$ primes $p,\  p\le
X.$ That is, we still have
$$
|\cE'|=\pi(X)(1+o(1)).
$$
Now,  the range of summation over $p$ in the above bound we concise
to $\cE'.$ Then
$$
\sum_{p\in \cE'}\max_{(a, p)=1}\Big|\sum_{n\le
T}\gamma_n\ep(a\lambda^{s_n})\Big|^2\ll \pi(X)T^2(\log
T)^{-\varepsilon/2}.
$$
From this, by taking $\gamma_n=1,$ we deduce the following
consequence.
\begin{cor}
\label{cor:Mersennetype} Let $s_n$ satisfy the condition
\eqref{eqn:s_n} and let $T=[X(\log X)^{2+\varepsilon}].$ Then the
inequality
$$
\max_{(a,p)=1}\Big|\sum_{n\le T}e^{2\pi
i\frac{a\lambda^{s_n}}{p}}\Big|\ll T(\log T)^{-\varepsilon/5}
$$
holds for all primes $p,\ p\le X,$ except at most $o(\pi(X))$ of
them.
\end{cor}

We recall that the {\it discrepancy} $D$ of a sequence of $N$ points
$(x_j)_{j=1}^N$ of the unit interval $[0,1)$ is  defined as
$$
D =\sup_{0\le a, b\le 1} \left|\frac{A(a,b) } {N}- (b-a)\right|,
$$
where $A(a, b)$ is the number of points of this sequence which
belong to $[a, b).$

Now let $D(p, X)$ denote the discrepancy of the fractional parts
$$
\{\lambda^{s_n}/p\},\quad n\le X(\log X)^{2+\varepsilon},
$$
where $s_n$ satisfies the condition~\eqref{eqn:s_n}. According to
the well-known Erd\H os-Tur{\'a}n relation between the discrepancy
and the associated exponential sums (see~\cite{DrTi}, or
alternatively one can use Theorem~4 of~\cite{Gar}), we derive from
Corollary~\ref{cor:Mersennetype} that for $\pi(X)(1+o(1))$ primes
$p, \ p\le X,$ the following bound holds with some
$\varepsilon_1>0:$
$$
D(p, X)\ll (\log X)^{-\varepsilon_1}.
$$
In other words, the numbers
$$
\lambda^{s_n},\quad n\le X(\log X)^{2+\varepsilon}
$$
are uniformly distributed modulo $p$ for any given $\varepsilon>0.$
In particular, one can take $s_n=[q_n^{c}],$ where $1\le c\le 15/14$
and $q_n$ denotes the $n$-th prime number.

The following Theorem is an analogy of Theorem~\ref{thm:main1},
where the range of summation over $n$ now depends on $p.$

\begin{theorem}
\label{thm:main3depend} Let $T_p, \ p\in \cE,$ be any positive
integers with $T_p\le T$ and let $\cE_1\subset \cE.$ For any
positive numbers $L$ and $K$  the following bound holds:
\begin{eqnarray*}
&&\sum_{p\in \cE_1}\frac{1}{\tau(p-1)}\max_{(a, p)=1}
\Big|\sum_{n\le T_p}\gamma_n\ep(a\lambda^{s_n})\Big|^2\\
&&\ll \left(X+s_TX^{1/7}\Delta^{-3/7}L+TL^{-7/4}\right)XT(\log
K)^{2}+\frac{T^2}{K^2}\sum_{p\in \cE_1}\frac{1}{\tau(p-1)}.
\end{eqnarray*}
\end{theorem}

Taking $\cE_1=\cE$ and $K=T$ and observing that the last term never
dominates, we see that Theorem~\ref{thm:main3depend} extends
Theorem~\ref{thm:main1} to more general sums at the cost of the
slight factor $(\log T)^2.$ In some applications one can further
relax this factor by special choices of $\cE_1$ and $K.$

\bigskip

One may want to have an explicit estimate for $|\overline{\cE}|,$
where
$$
\overline{\cE} =\{p:\ p\le X, p\not\in \cE\}.
$$
In this connection we remark that the argument given
in~\cite{ErdMur} immediately shows the inequality
$$
|\overline{\cE}|\ll\frac{\Delta^2}{\log \Delta}.
$$
Indeed
$$
\prod_{p\in\overline{\cE}}p \ | \prod_{k\le \Delta}(\lambda^{k}-1),
$$
Therefore, if $\omega(n)$ denotes the number of prime divisors of
$n,$ then we have
$$
|\overline{\cE}|\ll \omega\left(\prod_{k\le
\Delta}(\lambda^{k}-1)\right)\ll \frac{\Delta^2}{\log \Delta},
$$
where we have used the well known bound  $\omega(n) \ll (\log
n)(\log \log n)^{-1}.$

For $\Delta=X^{1/2+o(1)}$ one can use the results of~\cite{IndlTim}.

\section{Lemmas}

We need the version of the large sieve inequality applied to our
situation (recall that $|\gamma_n|\le 1$).

\begin{lemma}
\label{lem:largesieve} For any $K\ge 1$ the following estimate
holds:
$$
\sum_{k\le K}\sum_{\substack{1\le c\le k\\
(c,k)=1}}\left|\sum_{n\le T}\gamma_n\ek(cs_n)\right|^2\ll
(K^2+s_T)T.
$$
\end{lemma}

For the proof, see for example, \cite[pp. 153-154]{Da}.

\bigskip

We also recall the following bound of Heath-Brown and
Konyagin~\cite{HBK}.

\begin{lemma}
\label{lem:HeathKon} Let an integer $\theta$ be of multiplicative
order $t$ modulo $p$. Then the following  bound holds:
$$
\max_{(a, p) = 1} \left|\sum_{z=1}^t \ep(a \theta^z) \right| \ll
\min\{ p^{1/2},   p^{1/4}t^{3/8},   p^{1/8}t^{5/8}\}.
$$
\end{lemma}

Instead of Lemma~\ref{lem:HeathKon} one can use the bound due to
Bourgain-Konyagin~\cite{BourKon}, which however does not improve our
final results.

\section{Proof of Theorem~\ref{thm:main1}}

If $L\le 1,$ then the estimate of Theorem~\ref{thm:main1} becomes
trivial. Therefore, we will suppose that $L>1.$

Denote
$$
\sigma_p(a)=\sum_{n\le T}\gamma_n\ep(a\lambda^{s_n}).
$$
We have
$$
\sigma_p(a)=\sum_{x=1}^{t_p}\sum_{\substack{n\le T\\{t_p}|
s_n-x}}\gamma_n\ep(a\lambda^{s_n})=\frac{1}{t_p}\sum_{x=1}^{t_p}\sum_{b=1}^{t_p}\sum_{n\le
T}\gamma_n\et(b(s_n-x))\ep(a\lambda^x).
$$
For each divisor $d|t_p$ we collect together the values of $b$ with
$(b,t_p)=d.$ Thus
$$
\sigma_p(a)=\frac{1}{t_p}\sum_{d|t_p}\sum_{x=1}^{t_p}{\sum_{\substack{c\le t_p/d\\
(c,t_p/d)=1}}}\sum_{n\le T}\gamma_n\etd(c(s_n-x))\ep(a\lambda^x).
$$

We treat the cases of big and small values of $d$ separately. For
big values of $d$ we will enjoy the summation over $x$ in a proper
way to get sufficient to our purposes cancellations. The small
values of $d$ are treated in a different way. Thus, we define
$$
v_p=t_p^{4/7}p^{1/7}
$$
and set
\begin{equation}
\label{eqn:R1}
R_1=\max_{(a,p)=1}\left|\frac{1}{t_p}\sum_{\substack{d|t_p\\
d\ge Lv_p}}\sum_{x=1}^{t_p}{\sum_{\substack{c\le t_p/d\\
(c,t_p/d)=1}}}\sum_{n\le
T}\gamma_n\etd(c(s_n-x))\ep(a\lambda^x)\right|,
\end{equation}
\begin{equation}
\label{eqn:R2}
R_2=\max_{(a,p)=1}\left|\frac{1}{t_p}\sum_{\substack{d|t_p\\
d<Lv_p}}\sum_{x=1}^{t_p}{\sum_{\substack{c\le t_p/d\\
(c,t_p/d)=1}}}\sum_{n\le
T}\gamma_n\etd(c(s_n-x))\ep(a\lambda^x)\right|.
\end{equation}
Then
$$
\max_{(a,p)=1}|\sigma_p(a)|\le R_1+R_2.
$$
In particular,
\begin{equation}
\label{eqn:sigmasquaretau} \sum_{p\in
\cE}\frac{1}{\tau(p-1)}\max_{(a,p)=1}|\sigma_p(a)|^2\le \sum_{p\in
\cE}\frac{R_1^2}{\tau(p-1)}+\sum_{p\in \cE}\frac{R_2^2}{\tau(p-1)}.
\end{equation}

Our aim is to estimate the sums on the right hand side
of~\eqref{eqn:sigmasquaretau}.

To estimate $R_1$ we divide the interval of summation over $x$ to
progressions of the form $y+zt_p/d,\ 1\le y\le t_p/d,\ 1\le z\le d.$
Thus
$$
R_1=\max_{(a,p)=1}\left|\frac{1}{t_p}\sum_{\substack{d|t_p\\
d\ge Lv_p}}\sum_{y=1}^{t_p/d}\sum_{z=1}^{d}{\sum_{\substack{c\le t_p/d\\
(c,t_p/d)=1}}}\sum_{n\le
T}\gamma_n\etd(c(s_n-y))\ep(a\lambda^y\lambda^{zt_p/d})\right|,
$$
whence
$$
R_1\ll \frac{1}{t_p}\sum_{\substack{d|t_p\\
d\ge Lv_p}}\sum_{y=1}^{t_p/d}\left|{\sum_{\substack{c\le t_p/d\\
(c,t_p/d)=1}}}\sum_{n\le
T}\gamma_n\etd(c(s_n-y))\right|\max_{(a,p)=1}\left|\sum_{z=1}^{d}\ep(a\lambda^y\lambda^{zt_p/d})\right|.
$$
 The sum over $z$ is estimated by Lemma~\ref{lem:HeathKon}.
Since $\lambda^{t_p/d}$ is an element of multiplicative order $d,$
then from Lemma~\ref{lem:HeathKon} we derive
$$
\max_{(a,p)=1}\left|\sum_{z=1}^{d}\ep(a\lambda^y\lambda^{zt_p/d})\right|\ll
p^{1/8}d^{5/8}.
$$
Therefore,
\begin{equation}
\label{eqn:R1R3}
R_1\ll \sum_{\substack{d|t_p\\
d\ge Lv_p}}p^{1/8}d^{5/8}R_3,
\end{equation}
where
$$
R_3= \frac{1}{t_p}\sum_{y=1}^{t_p/d}\left|{\sum_{\substack{c\le t_p/d\\
(c,t_p/d)=1}}}\sum_{n\le T}\gamma_n\etd(c(s_n-y))\right|.
$$
Next, applying the Cauchy inequality we obtain
\begin{eqnarray*}
&&R_3^2\ll\frac{1}{dt_p}\sum_{y=1}^{t_p/d}\left|{\sum_{\substack{c\le t_p/d\\
(c,t_p/d)=1}}}\sum_{n\le T}\gamma_n\etd(c(s_n-y))\right|^2=\\
&&\frac{1}{dt_p}\sum_{y=1}^{t_p/d}\sum_{\substack{c_1\le t_p/d\\
(c_1,t_p/d)=1}}\sum_{\substack{c_2\le t_p/d\\
(c_2,t_p/d)=1}} \sum_{n_1\le T}\sum_{n_2\le
T}\gamma_{n_1}\overline{\gamma}_{n_2}\etd(c_1(s_{n_1}-y)-c_2(s_{n_2}-y)).
\end{eqnarray*}
Observe that
$$
\sum_{y=1}^{t_p/d} \etd(-c_1y+c_2y)= \left\{
\begin{array}{ll}
t_p/d,& \quad \mbox{if}\ c_1\equiv c_2 \pmod {t_p/d}, \\
0,& \quad \mbox{if}\ c_1 \not \equiv c_2 \pmod {t_p/d}.
\end{array} \right.
$$
Hence,
$$
R_3^2\ll\frac{1}{d^2}\sum_{\substack{c\le t_p/d\\
(c,t_p/d)=1}}\sum_{n_1\le T}\sum_{n_2\le
T}\gamma_{n_1}\overline{\gamma}_{n_2}\etd(c(s_{n_1}-s_{n_2})).
$$
Estimating trivially the sums over $c,\  n_1$ and $n_2$ we obtain
$$
R_3^2 \ll \frac{t_p}{d^3}T^2.
$$
Substituting this in~\eqref{eqn:R1R3}, we derive that
$$
R_1^2\ll \tau(p-1)T^2\sum_{\substack{d|t_p\\
d\ge Lv_p}} \frac{p^{1/4}t_p}{d^{7/4}}.
$$
Since $v_p=t_p^{4/7}p^{1/7},$ then
$$
R_1^2\ll
\tau(p-1)^2T^2\frac{p^{1/4}t_p}{(Lv_p)^{7/4}}=\tau(p-1)^2T^2L^{-7/4},
$$
whence
$$
\frac{R_1^2}{\tau(p-1)}\ll \tau(p-1)T^2L^{-7/4}.
$$
Application of the  Titchmarsh estimate~\eqref{eqn:titch} yields
\begin{equation}
\label{eqn:finR1tau} \sum_{p\in \cE}\frac{R_1^2}{\tau(p-1)}\ll
XT^2L^{-7/4}.
\end{equation}

Now we proceed to treat $R_2.$ From~\eqref{eqn:R2} we have
$$
R_2\le \frac{1}{t_p}\sum_{\substack{d|t_p\\
d<Lv_p}}\sum_{x=1}^{t_p}\left|{\sum_{\substack{c\le t_p/d\\
(c,t_p/d)=1}}}\sum_{n\le T}\gamma_n\etd(c(s_n-x))\right|.
$$
We apply the Cauchy inequality to the sums over $d$ and $x$ and then
obtain
$$
R_2^2\ll\frac{\tau(p-1)}{t_p}\sum_{\substack{d|t_p\\
d<Lv_p}}\sum_{x=1}^{t_p}\left|{\sum_{\substack{c\le t_p/d\\
(c,t_p/d)=1}}}\sum_{n\le T}\gamma_n\etd(c(s_n-x))\right|^2,
$$
whence
\begin{eqnarray*}
&&\frac{R_2^2}{\tau(p-1)}\ll\\
&&\frac{1}{t_p}\sum_{\substack{d|t_p\\
d<Lv_p}}\sum_{x=1}^{t_p}\sum_{\substack{c_1\le t_p/d\\
(c_1,t_p/d)=1}}\sum_{\substack{c_2\le t_p/d\\
(c_2,t_p/d)=1}} \sum_{n_1\le T}\sum_{n_2\le
T}\gamma_{n_1}\overline{\gamma}_{n_2}\etd(c_1(s_{n_1}-x)-c_2(s_{n_2}-x)).
\end{eqnarray*}
The summation over $x$ guarantees that $c_1=c_2.$ Therefore,
$$
\frac{R_2^2}{\tau(p-1)}\ll \sum_{\substack{d|t_p\\
d<Lv_p}}\sum_{\substack{c\le t_p/d\\
(c,t_p/d)=1}}\sum_{n_1\le T}\sum_{n_2\le
T}\gamma_{n_1}\overline{\gamma}_{n_2}\etd(c(s_{n_1}-s_{n_2})),
$$
whence
$$
\frac{R_2^2}{\tau(p-1)}\ll \sum_{\substack{d|t_p\\
d<Lv_p}}\sum_{\substack{c\le t_p/d\\
(c,t_p/d)=1}}\left|\sum_{n\le T}\gamma_n\etd(cs_n)\right|^2.
$$
Summing up both sides of this bound  over $p\in \cE,$ we obtain
$$
\sum_{p\in \cE}\frac{R_2^2}{\tau(p-1)}\ll \sum_{p\in \cE}\sum_{\substack{d|t_p\\
d<Lv_p}}\sum_{\substack{c\le t_p/d\\
(c,t_p/d)=1}}\left|\sum_{n\le T}\gamma_n\etd(cs_n)\right|^2.
$$

We divide the interval $(\Delta, X]$ into disjoint subintervals
$(X_j, X_{j+1}],$ where
$$
X_1=\Delta,\quad X_{j+1}=\min \{2X_j,X\}.
$$
Denote by $\cE_j$ the subset of $\cE$ such that $t_p\in (X_j,
X_{j+1}]$ for any $p\in \cE_j.$ Next, define
$$
V_j=2X_j^{4/7}X^{1/7}
$$
and observe that $V_j$ does not depend on $p,$ and $V_j\ge v_p$ for
any $p\in \cE_j.$ Thus,
$$
\sum_{p\in \cE}\frac{R_2^2}{\tau(p-1)}\ll \sum_{j}\sum_{p\in
\cE_j}\sum_{\substack {d|t_p\\ d<LV_j}}\sum_{\substack{c\le t_p/d \\
(c, t_p/d)=1}}\Big|\sum_{n\le T}\gamma_n\etd(cs_n)\Big|^2.
$$
We remember that  $j\ll \log X,$ $2^jX_1\ll X$ and
$$
\Delta\le X_j<X_{j+1}\le 2X_j\le 2X.
$$
Observe that for different primes $p, p\in \cE_j,$ the corresponding
values of $t_p$ do not have to be different. For a given $r\in (X_j,
X_{j+1}]$ denote by $s(r)$ the number of all primes $p, p\in \cE_j,$
for which $t_p=r.$ Since $ p-1\equiv 0\pmod r,$ then
$$
s(r)\le X/r\le X/X_j.
$$
Therefore,
$$
\sum_{p\in \cE}\frac{R_2^2}{\tau(p-1)}\ll
\sum_j\frac{X}{X_j}\sum_{r\in(X_j, X_{j+1}]}
\sum_{\substack {d|r\\ d<LV_j}}\sum_{\substack{c\le r/d \\
(c, r/d)=1}}\Big|\sum_{n\le T}\gamma_n\erd(cs_n)\Big|^2.
$$
Changing the order of summation over $r$ and $d$ we deduce
\begin{equation}
\label{eqn:Rd} \sum_{p\in \cE}\frac{R_2^2}{\tau(p-1)}\ll\sum_j
\frac{X}{X_j}\sum_{d< LV_j}F_j(d),
\end{equation}
where
\begin{eqnarray*}
&&F_j(d)=\sum_{\substack{r\in (X_j,
X_{j+1}]\\ r\equiv 0\pmod d}} \sum_{\substack{1\le c\le r/d \\
(c, r/d)=1}}\Big|\sum_{n\le T}\gamma_n\erd(cs_n)\Big|^2 \\
&&=\sum_{k\in (X_jd^{-1},
X_{j+1}d^{-1}]} \, \sum_{\substack{1\le c\le k \\
(c, k)=1}}\Big|\sum_{n\le T}\gamma_n\ek(cs_n)\Big|^2.
\end{eqnarray*}

To estimate $F_j(d)$ we apply the large sieve inequality given in
Lemma~\ref{lem:largesieve}. Then
\begin{eqnarray*}
F_j(d) \ll (X_j^2d^{-2}+s_T)T.
\end{eqnarray*}
Inserting this bound into~\eqref{eqn:Rd}, we obtain
$$
\sum_{p\in \cE}\frac{R_2^2}{\tau(p-1)}\ll
\sum_j\frac{X}{X_j}\sum_{d< LV_j}(X_j^2d^{-2}+s_T)T,
$$
whence
$$
\sum_{p\in \cE}\frac{R_2^2}{\tau(p-1)}\ll
\sum_jX(X_j+s_TV_jLX_j^{-1})T.
$$
Since $ V_j=2 X_j^{4/7}X^{1/7},$ we have
$$
\sum_{p\in \cE}\frac{R_2^2}{\tau(p-1)}\ll
XT\left(\sum_jX_j+s_TLX^{1/7}\sum_j X_j^{-3/7}\right).
$$
Finally, from the definition of $X_j$ we know that
$$
\sum_jX_j\ll X,\qquad \sum_{j}X_j^{-3/7}\ll \Delta^{-3/7}.
$$
Therefore,
\begin{equation}
\label{eqn:finR2} \sum_{p\in \cE}\frac{R_2^2}{\tau(p-1)}  \ll
XT(X+s_TX^{1/7}\Delta^{-3/7}L).
\end{equation}

Theorem~\ref{thm:main1} now follows from \eqref{eqn:sigmasquaretau},
\eqref{eqn:finR1tau} and \eqref{eqn:finR2}.

\section{Proof of Theorem~\ref{thm:main3depend}}

For $K\le 10$ the estimate of Theorem~\ref{thm:main3depend} is
trivial. Therefore, we will suppose that $K>10.$

Set $M=[T/K].$ Without loss of generality we may assume that for
$n\ge 1,$
$$
\gamma_{T+n}=0, \quad s_{T+n}=s_T+n.
$$
Applying the shifting argument we obtain
\begin{equation}
\label{eqn:shift}
\left|\sum_{n=1}^{T_p}\gamma_n\ep(a\lambda^{s_n})\right|^2\ll
\frac{1}{(M+1)^2}\left|\sum_{n=1}^{T_p}\sum_{r=0}^{M}\gamma_{n+r}\ep(a\lambda^{s_{n+r}})\right|^2+
\frac{T^2}{K^2}.
\end{equation}
Further, we have
\begin{equation}
\label{eqn:smooth}
\sum_{n=1}^{T_p}\sum_{r=0}^{M}\gamma_{n+r}\ep(a\lambda^{s_{n+r}})=\frac{1}{2
T+1}\sum_{b=-T}^{T}\sum_{m=1}^{2T}\sum_{r=0}^{M}\sum_{n=1}^{T_p}\gamma_me^{2\pi
i\frac{b(n+r-m)}{2T+1}}\ep(a\lambda^{s_m}).
\end{equation}
By the Cauchy inequality,
\begin{eqnarray*}
&& \left(\sum_{0<|b|\le
T}\left|\sum_{n=1}^{T_p}\sum_{r=0}^{M}e^{2\pi
i\frac{b(n+r)}{2T+1}}\right|\left|\sum_{m=1}^{2T}\gamma_me^{2\pi
i\frac{bm}{2T+1}}\ep(a\lambda^{s_m})\right|\right)^2\ll
\\
&& \left(\sum_{0<|b|\le
T}\left|\sum_{n=1}^{T_p}\sum_{r=0}^{M}e^{2\pi
i\frac{b(n+r)}{2T+1}}\right|\right)\times \\
&& \left(\sum_{0<|b|\le
T}\left|\sum_{n=1}^{T_p}\sum_{r=0}^{M}e^{2\pi
i\frac{b(n+r)}{2T+1}}\right|\left|\sum_{m=1}^{T}\gamma_me^{2\pi
i\frac{bm}{2T+1}}\ep(a\lambda^{s_m})\right|^2\right).
\end{eqnarray*}
Hence, using
$$
\left|\sum_{n=1}^{T_p}e^{2\pi i\frac{bn}{2T+1}}\right|\ll
\frac{T}{|b|},
$$
we obtain the bound
\begin{eqnarray*}
\left(\sum_{0<|b|\le T}\left|\sum_{n=1}^{T_p}\sum_{r=0}^{M}e^{2\pi
i\frac{b(n+r)}{2T+1}}\right|\left|\sum_{m=1}^{2T}\gamma_me^{2\pi
i\frac{bm}{2T+1}}\ep(a\lambda^{s_m})\right|\right)^2\ll
\\
T^2\left(\sum_{0<|b|\le
T}\frac{|S(b)|}{|b|}\right)\left(\sum_{b=1}^{T}\frac{|S(b)|}{|b|}\left|\sum_{m=1}^{T}\gamma_me^{2\pi
i\frac{bm}{2T+1}}\ep(a\lambda^{s_m})\right|^2\right),
\end{eqnarray*}
where
\begin{equation}
\label{eqn:defS(b)} S(b)=\sum_{r=0}^{M}e^{2\pi i\frac{br}{2T+1}}.
\end{equation}
Combining this with~\eqref{eqn:shift} and~\eqref{eqn:smooth}, we
deduce
\begin{eqnarray*}
&& \left|\sum_{n=1}^{T_p}\gamma_n\ep(a\lambda^{s_n})\right|^2 \ll
\\ && \frac{1}{(M+1)^2}\left(\sum_{0<|b|\le
T}\frac{|S(b)|}{|b|}\right)\left(\sum_{0<|b|\le
T}\frac{|S(b)|}{|b|}\left|\sum_{m=1}^{T}\gamma_me^{2\pi
i\frac{bm}{2T+1}}\ep(a\lambda^{s_m})\right|^2\right)\\ &&+
\left|\sum_{m=1}^{T}\gamma_m\ep(a\lambda^{s_m})\right|^2+\frac{T^2}{K^2}
\end{eqnarray*}
Now we take maximum over $a,\ (a,p)=1,$ and observe that the maximum
of sums is not greater than the sum of maximums. We then divide the
estimate by $\tau(p-1)$ and perform the summation over $p\in\cE_1.$
This yields
\begin{eqnarray*}
&&\sum_{p\in
\cE_1}\frac{1}{\tau(p-1)}\max_{(a,p)=1}\left|\sum_{n=1}^{T_p}\gamma_n\ep(a\lambda^{s_n})\right|^2 \ll\\
&& \frac{1}{(M+1)^2}\left(\sum_{0<|b|\le
T}\frac{|S(b)|}{|b|}\right)\times \\ && \left(\sum_{0<|b|\le
T}\frac{|S(b)|}{|b|}
\sum_{p\in\cE_1}\frac{1}{\tau(p-1)}\max_{(a,p)=1}\left|\sum_{m=1}^{T}\gamma_me^{2\pi
i\frac{bm}{2T+1}}\ep(a\lambda^{s_m})\right|^2\right) \\
&&
+\sum_{p\in\cE_1}\frac{1}{\tau(p-1)}\max_{(a,p)=1}\left|\sum_{m=1}^{T}\gamma_m\ep(a\lambda^{s_m})\right|^2+
\frac{T^2}{K^2}\sum_{p\in\cE_1}\frac{1}{\tau(p-1)}.
\end{eqnarray*}
For each $b$ to the sum
$$
\sum_{p\in\cE_1}\frac{1}{\tau(p-1)}\max_{(a,p)=1}\left|\sum_{m=1}^{T}\gamma_me^{2\pi
i\frac{bm}{2T+1}}\ep(a\lambda^{s_m})\right|^2
$$
we apply Theorem~\ref{thm:main1} with  $\gamma_n$ substituted by
$\gamma_ne^{2\pi i\frac{bn}{2T+1}}.$ Thus,
\begin{eqnarray*}
&&\sum_{p\in
\cE_1}\frac{1}{\tau(p-1)}\max_{(a,p)=1}\left|\sum_{n=1}^{T_p}\gamma_n\ep(a\lambda^{s_n})\right|^2 \ll\\
&&
\left(\frac{1}{(M+1)^2}\left(\sum_{b=1}^{T}\frac{|S(b)|}{b}\right)^2+
1\right)\left(X+s_TX^{1/7}\Delta^{-3/7}L+TL^{-7/4}\right)XT\\
&& + \frac{T^2}{K^2}\sum_{p\in\cE_1}\frac{1}{\tau(p-1)}.
\end{eqnarray*}
Now it remains to prove that
$$
\sum_{b=1}^{T}\frac{|S(b)|}{b}\ll (M+1)\log K.
$$
To this end, choose $\ell=[\log K]$ and use the Holder inequality to
obtain
\begin{equation}
\label{eqn:holderS(b)/b} \sum_{b=1}^{T}\frac{|S(b)|}{b}\le
\left(\sum_{b=1}^{2T+1}\frac{1}{b^{2\ell/(2\ell-1)}}\right)^{1-1/2\ell}
\left(\sum_{b=1}^{2T+1}|S(b)|^{2\ell}\right)^{1/2\ell}.
\end{equation}
Next, we have
\begin{equation}
\label{eqn:estsumb}
\sum_{b=1}^{2T+1}\frac{1}{b^{2\ell/(2\ell-1)}}\ll
\int_{1}^{\infty}x^{-1-(2\ell-1)^{-1}}dx = 2\ell-1 \ll \log K.
\end{equation}
Besides, from the definition of $S(b),$ see~\eqref{eqn:defS(b)}, it
follows
\begin{equation}
\label{eqn:S(b)tocongr} \sum_{b=1}^{2T+1}|S(b)|^{2\ell}=(2T+1)J,
\end{equation}
where $J$ denotes the number of solutions to the congruence
$$
\sum_{i=1}^{\ell}x_i\equiv\sum_{i=1}^{\ell}y_i\pmod{(2T+1)},\qquad
0\le x_i, y_j\le M.
$$
Since $M<T,$ then the trivial estimate gives $J\le (M+1)^{2\ell-1}.$
Besides, $T<K(M+1).$ Therefore,
$$
\sum_{b=1}^{2T+1}|S(b)|^{2\ell}\ll K(M+1)^{2\ell},
$$
whence, in view
of~\eqref{eqn:holderS(b)/b}--\eqref{eqn:S(b)tocongr}, we conclude
that
$$
\sum_{b=1}^{T}\frac{|S(b)|}{b}\ll (\log K) (M+1)K^{1/(2\ell)}\ll
(M+1)\log K.
$$

\section{Exponent pairs for Gauss sums}

We remark that if in Lemma~\ref{lem:HeathKon} we have the bound
\begin{equation}
\label{eqn:pairs} \max_{(a, p) = 1} \left|\sum_{z=1}^t \ep(a
\theta^z) \right| \ll p^{\alpha}t^{\beta}
\end{equation}
with $0\le\alpha, \beta\le 1,$ then the right hand side of the
estimate of Theorem~\ref{thm:main1} can be substituted by
$$
(X+s_TX^{\frac{2\alpha}{3-2\beta}}\Delta^{-\frac{2-2\beta}{3-2\beta}}L+TL^{-3+2\beta})XT.
$$
In particular, Corollary~\ref{cor:Mersennetype} takes place for the
sequence $s_n$ satisfying the condition
$$
s_T\le T^{1+\frac{1-2\alpha-\beta}{3-2\beta}+o(1)}.
$$
Define $\cK$ to be the set of all ordered pairs $\{\alpha,\beta\}$
with $0\le \alpha, \beta\le 1$ and satisfying the
property~\eqref{eqn:pairs}. Konyagin~\cite{Kon} has proved that the
set $\cK$ contains the pair $\{\alpha_n, \beta_n\}$ defined as
$$
\alpha_n=\frac{1}{2n^2},\quad
\beta_n=1-\frac{2}{n^2}+\frac{1}{2^{n-1}n^2}
$$
for any positive integer $n.$ Furthermore, $\cK$ also contains the
pair $\{\alpha_n', \beta_n'\}$ given by
$$
\alpha_n'=\frac{1}{2n(n+1)},\quad
\beta_n'=1-\frac{2}{n(n+1)}+\frac{3}{2^{n+1}n(n+1)}.
$$

We now define the function $f:\cK\to \mathbb{R}$ by
$$
f(x,y)=1+\frac{1-2x-y}{3-2y}.
$$
The problem is to find the value of $f(x,y)$ as big as possible. The
result of the present paper corresponds to the pair $\{\alpha_2,
\beta_2\}$ (which is due to Heath-Brown and Konyagin). Other pairs
give less precise bounds. Next, we note that $\cK$ is a convex set.
That is, if
$$
\{\alpha, \beta\}\in \cK, \quad \{\alpha', \beta'\}\in \cK,
$$
then for any $x,$ $0\le x\le 1,$
$$
\{x\alpha+(1-x)\alpha', x\beta+(1-x)\beta'\}\in \cK.
$$
However, this property applied to any two given pairs, in particular
to the pairs due to Konyagin, is not sufficient to get further
improvements, and it would be very interesting, similar to the set
of exponent pairs, have more nontrivial properties of $\cK.$ The
truth of the conjecture of Montgomery, Vaughan and
Wooley~\cite{Mont} would imply
$$
\{\varepsilon, 1/2+\varepsilon\}\in \cK,
$$
which can be considered as an analogy of the exponent pair
hypothesis.

Finally, we remark that the method we have applied leads to the
following generalization of our main result.

\begin{theorem}
\label{thm:general} For any positive integer $N,$ any $L>0,$ any
pair $\{\alpha, \beta\}\in\cK$ and any complex coefficients
$\delta_n, \, 1\le n\le N,$ the following bound holds:
\begin{eqnarray*}
&&\sum_{p\in
\cE}\frac{1}{\tau(p-1)}\max_{(a,p)=1}\left|\sum_{n=1}^N\delta_n\ep(a\lambda^n)\right|^2\ll
\\&&
X(X+NX^{\frac{2\alpha}{3-2\beta}}\Delta^{-\frac{2-2\beta}{3-2\beta}}L)
\sum_{n=1}^N|\delta_n|^2+
XL^{-3+2\beta}\left(\sum_{n=1}^N|\delta_n|\right)^2,
\end{eqnarray*}
where the implied constant depends only on the pair
$\{\alpha,\beta\}.$
\end{theorem}

\end{document}